# Boundary Element Method for the Dirichlet Problem for Laplace's Equation on a Disk

**Misael M. Morales**[1*†]
**Shirley Pomeranz**[1]
1. The University of Tulsa, Tulsa, OK, USA.
† *now at The University of Texas at Austin, Austin, TX, USA.*
\* *corresponding author: misaelmorales@utexas.edu*

The Boundary Element Method (BEM) is implemented using piecewise linear elements to solve the two-dimensional Dirichlet problem for Laplace's equation posed on a disk. A benefit of the BEM as opposed to many other numerical solution techniques is that discretization only occurs on the boundary, i.e., the complete domain does not need to be discretized. This provides an advantage in terms of time and cost. The algorithm's performance is illustrated through sample test problems with known solutions. A comparison between the exact solution and the BEM numerical solution is done, and error analysis is performed on the results.

## ■ Introduction

Our study of the Boundary Element Method (BEM) will focus on solving the Dirichlet problem for Laplace's equation posed on a unit disk centered at the origin. This is a common partial differential equation (PDE) that arises in several fields of science and engineering, including fluid flow in porous media, geomechanics, electromagnetics, thermodynamics, and more [1, 2]. It is especially good for methods with infinite or semi-infinite domains, where discretization at the boundary is less expensive than the full-domain discretization required by other methods such as the finite element method and the finite dif-ference method. Like its counterparts, the BEM can also be applied to other types of PDEs but requires special treatment of time derivatives including Laplace transforms and numerical inversion [3, 4, 5].

For our 2D elliptic problem, the model problem is given by

$$\begin{aligned}\nabla^2 u = \partial_{x,x} u + \partial_{y,y} u = 0, &\quad (x, y) \in \Omega, \\ u = f(x, y), &\quad (x, y) \in \partial\Omega,\end{aligned} \quad (1)$$

where, for our problem, $\Omega$ is a closed, bounded domain with piecewise smooth boundary, $\partial\Omega$.

We will focus on solving the Dirichlet problem for Laplace's equation on a disk domain, taking advantage of the BEM. A Dirichlet boundary condition is required to find the solution; thus the value that a solution takes along the boundary is specified. Collocation will be used to obtain a system of equations. We require the fundamental solution, i.e., Freespace Green's Function [6], to Laplace's equation, given by

$$w(x, y) = \frac{-1}{2\pi} \ln\left(\sqrt{(x - \xi)^2 + (y - \eta)^2}\right), \quad (2)$$

for $(x, y) \neq (\xi, \eta)$, where $(x, y)$ is an arbitrary field point on the boundary and $(\xi, \eta)$ is an arbitrary source (load, collocation, evaluation) point in the domain or on the boundary. A disadvantage of the BEM is that a fundamental solution for the differential operator is required. Fundamental solutions are known for many basic differential operators; however, if a fundamental solution is not known, then the BEM cannot be used.

This fundamental solution is the analytical solution of Laplace's equation under the action of a point source on an infinite domain [1, 7]. For many elliptic equations, including Laplace's equation, fundamental solutions are well known in the literature [8]. The fundamental solution serves as a weight factor in the integral equation formulation. By transforming our PDE into an integral equation, we formulate the foundations for the BEM numerical solution technique.

To implement the BEM, we multiply the PDE by the fundamental solution, here denoted $w$, and integrate over the domain. The solution for the unknown boundary outward normal flux, denoted $\frac{\partial u}{\partial n}$, and, subsequently the interior primary unknown, denoted $u$, is obtained using Green's identities (i.e., integrate by parts twice and apply the divergence theorem) [9],

$$\begin{aligned}\left(\nabla^2 u\right) w = 0 \Rightarrow \\ \int_\Omega \left(\nabla^2 u\right) w \, d\Omega = \int_\Omega \nabla \cdot ((\nabla u) w) - \nabla u \cdot \nabla w \, d\Omega = \\ \int_{\partial\Omega} \frac{\partial u}{\partial n} w \, ds - \int_{\partial\Omega} u \frac{\partial w}{\partial n} \, ds + \int_\Omega u \nabla^2 w \, d\Omega = 0,\end{aligned} \quad (3)$$

where $ds$ is the differential element of length with respect to $(x, y)$ along the the bound-ary, $\partial\Omega$.

We now have the basic boundary integral equation. This first provides the BEM equation that is solved for the unknown outward normal boundary flux at each boundary collocation point, and this equation is given by

$$c(\xi, \eta)\, u(\xi, \eta) + \int_{\partial\Omega} u(x, y)\, \frac{\partial w(x, y, \xi, \eta)}{\partial n}\, ds = \int_{\partial\Omega} w(x, y, \xi, \eta)\, \frac{\partial u(x, y)}{\partial n}\, ds, \quad (\xi, \eta) \in \partial\Omega, \tag{4}$$

where the flux operator, $\frac{\partial}{\partial n}$, denotes the outward normal partial derivative taken with respect to the field coordinates. Here the factor $c(\xi, \eta)$ is known as the free-term coefficient, and the collocation points lie on the boundary. The primary unknown, $u$, and its normal derivative, $\frac{\partial u}{\partial n}$, are to be discretized. The integration is performed with respect to the field point coordinates, $(x, y)$. It is also important to note that even though the discretization of the domain is only occurring on the boundary, we still require evaluation of the primary unknown in the interior of the domain.

To find the solution for the primary unknown at points in the domain, select arbitrary interior collocation points. Rearrange Eq. (4) with, in this case, $c(\xi, \eta) = 1$ [10],

$$u(\xi, \eta) = -\int_{\partial\Omega} u(x, y)\, \frac{\partial w(x, y, \xi, \eta)}{\partial n}\, ds + \int_{\partial\Omega} \frac{\partial u(x, y)}{\partial n}\, w(x, y, \xi, \eta)\, ds, \tag{5}$$

$(\xi, \eta) \in \Omega$.

Eqs. (4) and (5) are discretized. The discretized boundary element system of equations can be expressed in either equation or matrix-vector form. We choose to use equation form in our code so that the code is more easily read. However, the matrix-vector form is more efficient in Mathematica.

For example, discretizing Eq. (5) gives

$$u_i = -\sum_{j=1}^{N} h_{i,j}\, u_j + \sum_{j=1}^{N} g_{i,j}\, q_j, \tag{6}$$

where $N$ is a positive integer denoting the number of boundary elements, $u_j$ is the known primary quantity evaluated at the $j^{\text{th}}$ boundary node, and $q_j = \left(\frac{\partial u}{\partial n}\right)_j$ is the unknown flux evaluated at the $j^{\text{th}}$ boundary node. The approximating integrals are

$$h_{i,j} = \int_{\Gamma_{j-1}} \frac{\partial w(x, y, \xi_i, \eta_i)}{\partial n}\, \beta_2^{j-1}(x, y)\, d\Gamma + \int_{\Gamma_j} \frac{\partial w(x, y, \xi_i, \eta_i)}{\partial n}\, \beta_1^j(x, y)\, d\Gamma \tag{7}$$

$$g_{i,j} = \int_{\Gamma_{j-1}} w(x, y, \xi_i, \eta_i)\, \beta_2^{j-1}(x, y)\, d\Gamma + \int_{\Gamma_j} w(x, y, \xi_i, \eta_i)\, \beta_1^j(x, y)\, d\Gamma, \tag{8}$$

where $\Gamma_j$ is the $j^{th}$ boundary element and $\beta_1^j(x, y_1)$ and $\beta_2^{j-1}(x, y)$ are the two non-zero basis functions at the $j^{th}$ boundary node. The superscripts for these basis functions in Eqs. (7) and (8) refer to the two consecutive boundary elements, indexed $j-1$ and $j$, which share the $j^{th}$ boundary node. The subscripts $i$ and $j$ in Eqs. (7) and (8) refer respectively to the $i^{th}$ collocation node and the $j^{th}$ boundary node, where the $j^{th}$ boundary element which is delineated by boundary nodes $j$ and $j+1$ [10]. As usual, the integrations are performed with respect to the $(x, y)$ coordinate.

The indexing in Eqs. (6), (7), and (8) uses a node-by-node assembly process on the boundary. Alternatively, assembly on the boundary could be organized element-by-element, as is done in the following code.

We obtain a dense matrix-vector problem, which can be computationally expensive to solve if the number of collocation points is large, but with the advantage that we end up solving a BEM problem with one dimension less than the original PDE formulation of the problem. Instead of discretizing and meshing in the two-dimensional space for the entire domain, we only need to discretize and mesh for one-dimensional integrations along the boundary. We can use standard numerical integration methods, preprocessors, and other strategies to find numerical solutions of our original PDEs. Even though we will create a grid in the domain interior to obtain evaluation points for the primary unknown, this is done for convenience here. It is not necessary to form a grid in the interior; instead, interior evaluation points for the primary unknown can be selected arbitrarily, and one BEM equation simply needs to be evaluated for each evaluation point.

In this notebook, we use Wolfram Mathematica - Student Edition, version 12.0.0.0.

# ■ BEM Implementation

## ▫ Notebook Introduction

### ▪ *Problem Definition*

We will solve a Dirichlet problem for Laplace's Equation on a disk bounded by the unit circle centered at the origin.

$$\nabla^2 u = \frac{\partial^2 u}{\partial x^2} + \frac{\partial^2 u}{\partial y^2} = \partial_{x,x} u + \partial_{y,y} u = u_{xx} + u_{yy} = 0, \qquad (x, y) \in \Omega, \tag{9}$$

$$u = f(x, y), \qquad (x, y) \in \partial\Omega.$$

For our problem, we will use uniform, piecewise linear boundary elements to approximate the circular boundary, the Dirichlet boundary conditions, and the outward normal boundary flux. These basis functions are sometimes called the "shape functions" [10].

This notebook serves as a template for solving related problems in science and engineering. The overall technique of the BEM is emphasized, but this notebook structure can be adapted as a template for a variety of problems. Some command outputs have been suppressed to make the code more readable and to save space.

### ▪ *Notebook Initialization*

To start the BEM template, the user should remove all global variables previously stored or start with a new Mathematica session. We save the current time is the variable *StartTime* in order to compute the overall notebook run time at the end of all computations.

```
In[1]:=  (* Remove["Global`*"] *)
```

```
In[2]:=  StartTime = AbsoluteTime[];
```

## ▫ Exact Solution

### ▪ *Problem Selection and Exact Solution*

In order check the accuracy of our BEM results, we use a test problem with a known exact analytical solution. We define five possible exact solutions, *uExact1*, *uExact2*, *uExact3*, *uExact4*, and *uExact5*. For example, the user can select one of the five exact solutions in or der to evaluate this notebook; *uExact1* is selected below.

```
In[3]:=  uExact1[{x_,y_}] = 1 + x² - y²;
         uExact2[{x_,y_}] = e^y Cos[x];
         uExact3[{x_,y_}] = 1 + Sin[π*x]*Sinh[π*y];
         uExact4[{x_,y_}] = 1 + -3*Cosh[4π]*Sin[2π*x]*Sinh[2π(y-2)] +
                            Csch[6π]*Sin[3π*x]*Sinh[3π*y];
         uExact5[{x_,y_}] = π*e^y*Cos[x - π/7] + e^(1-π*x) Cos[π*y - π/2] +
                            e^(5*x)/100 *Cos[5*y - π/2];
```

We rename the chosen exact solution as *uExact*. This will serve to demonstrate the BEM code in this notebook.

In[8]:= `uExact[{x_,y_}] = uExact1[{x,y}]`

Out[8]= $1 + x^2 - y^2$

To verify that the exact solution chosen from above satisfies Laplace's equation, we perform the following check, and ensure that $\partial_{x,x} + \partial_{y,y} = 0$.

In[9]:= `∂_{x,x}uExact[{x,y}] + ∂_{y,y}uExact[{x,y}]==0`

Out[9]= True

The exact outward normal boundary flux can be computed for this test problem, and we save it in *qExact*.

In[10]:= `qExact[{x_,y_}] = {∂_x uExact[{x,y}],∂_y uExact[{x,y}]}.{x,y}`

Out[10]= $2x^2 - 2y^2$

### ■ View the Exact Solution

To view the exact solution, we create a 3D plot in the problem domain, the disk with radius one centered at the origin, $x^2 + y^2 \le 1$.

```
In[11]:=  plotCircExact = Plot3D[uExact[{x,y}],{x,-1,1},{y,-1,1},
          PlotRange→All, ColorFunction→"ThermometerColors",
          ImageSize→Small,AxesLabel→{"x","y","u"},
          PlotLabel→"Exact Solution\non Domain Disk",
          RegionFunction→Function[{x,y,z},0<x^2+y^2≤1]]
```

Out[11]=

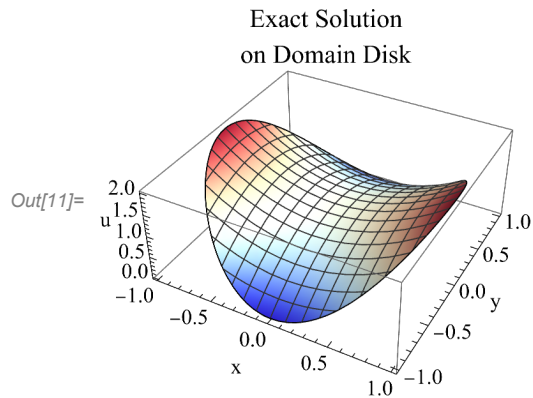

▲ **Figure 1.** The exact solution to the test problem in the problem domain, the unit disk centered at the origin.

## ◻ BEM: Fundamental Quantities and Definitions

### ▪ *Initialization Parameters*

The number of boundary nodes and interior nodes are specified by the user. For the circular boundary, we specify the number of uniformly spaced boundary nodes, *nBoundaryNodes*. For the interior of the disk, we specify the number of nodes that will create a Cartesian grid for the interior. For convenience in coding, the quantity *tempInteriorNodes* will be used later in the notebook to denote the number of nodes in the x-direction, with the same number of nodes in the y-direction. These nodes lie in the square that circumscribes the domain boundary circle. Therefore, the total number of nodes in the square Cartesian grid will be equal to the square of *tempInteriorNodes*. We will then select those inside the unit disk.

Even though we create a grid in the domain interior to obtain evaluation points for the primary unknown, this is done purely for convenience. It is not necessary to form a grid in the interior. Instead, interior evaluation points for the primary unknown can be selected arbitrarily, and one BEM equation simply needs to be evaluated for each evaluation point. Using a square geometry for the interior simplifies the generation of node coordinates and is done just for convenience so as to continue to work with rectangular instead of polar co-ordinates. That is, we chose to temporarily code on a square region and then only select collocation nodes inside the domain disk.

*In[12]:=*
```
nBoundaryNodes = 30;
tempInteriorNodes = 11;
```

The user can also specify the numerical integration technique desired to solve the problem. If left blank, Mathematica will implement the default strategy method for *NIntegrate*, which is a Global Adaptive scheme [11]. Several methods have been tested for our algorithm, including strategies, rules, and preprocessors. The strategies include "GlobalAdaptive", "LocalAdaptive", "DoubleExponential", "MultiPeriodic"; the quadrature rules include "TrapezoidalRule", "NewtonCotesRule", "GaussBerntsenEspelidRule", "GaussKronrodRule", "LobattoKronrodRule", "ClenshawCurtisRule", "MultiPanelRule", "CartesianRule", "MultiDimensionalRule"; the tested matrix preprocessors include "SymbolicPreprocessing", "SymbolicPiecewiseSubdivision", "EvenOddSubdivision", "UnitCubeRescaling".

Correct formatting of the numerical integration method is as follows

```
Method → {"preprocessor", Method → "rule"}
```

or

```
Method → {"strategy", Method → "rule"}
```

In this template we have selected the "SymbolicPiecewiseSubdivision" matrix preprocessor and the "NewtonCotesRule" quadrature rule. Newton-Cotes formulas are simple and efficient quadrature formulas, where the integrand is evaluated at equally spaced points, and include the common Trapezoidal and Simpson's rules. Specifying this method helps for speedy computations of the integrands. The preprocessor chosen, "Symbolic Piecewise Subdivision," partitions the integral into disjoint integration regions on each of which the integrand is continuous. Preprocessing the integrals and applying simple and efficient quadrature rules helps improve efficiency in timing and precision. Other methods can be used as well.

*In[14]:=*
```
integralMethod =
    {"SymbolicPiecewiseSubdivision",Method→"NewtonCotesRule"};
```

### ▪ Free Term Coefficient and Domain Variables

For this problem, the free term coefficient is given by [1, 6, 7]. The free term coefficient value is one when the collocation points are in the interior [10]. When the collocation points are on the boundary, for this problem and choice of uniform boundary element discretization, we have

$$c(\xi) = \begin{cases} \dfrac{(n-2)*\frac{\pi}{n}}{2\pi} = \dfrac{n-2}{2*n}, & \xi \in \partial\Omega, \\ 1 & \xi \in \Omega, \\ 0, & \xi \notin \Omega, \end{cases}$$

where *n* is the number of boundary nodes.

*In[15]:=*
```
FreeTermCoefficient = (nBoundaryNodes-2)/(2*nBoundaryNodes);
```

We define the quantities that will be used for the independent variables on the boundary elements, *xVar* and *yVar*. These quantities will be parameterized on each boundary element using the standardized variable $t \in [-1, 1]$. The $i^{\text{th}}$ boundary element is delineated by the nodes $i$ and $i+1$. We make sure to specify that the $(n + 1)^{\text{st}}$ node is the same as the $1^{\text{st}}$ node.

*In[16]:=*
```
Do[
 xVar[i,t_] = t(x[i+1]-x[i])/2 + (x[i]+x[i+1])/2;
 yVar[i,t_] = t(y[i+1]-y[i])/2 + (y[i]+y[i+1])/2,
 {i,1,nBoundaryNodes}]
```

*In[17]:=*
```
x[nBoundaryNodes+1] = x[1];
y[nBoundaryNodes+1] = y[1];
```

### ■ Basis Functions

The basis functions, also known as shape functions, are chosen as piecewise linear functions. Both the boundary element basis functions and fundamental quantities are parameterized using the parameter $t \in [-1, 1]$. We plot the basis functions and their derivatives.

*In[19]:=*
```
β1[t_] = (1-t)/2;

β2[t_] = (1+t)/2;
```

*In[21]:=*
```
plotBasisFunction1 =
Plot[{β1[t],β1'[t]},{t,-1,1},ImageSize→Small,
PlotStyle→{{Thickness[0.01],Red},
{Thickness[0.01],LightRed, Dashed}}];

plotBasisFunction2 =
Plot[{β2[t],β2'[t]},{t,-1,1},ImageSize→Small,
PlotStyle→{{Thickness[0.01],Blue},
{Thickness[0.01],LightBlue,Dashed}}];

Show[plotBasisFunction1,plotBasisFunction2,Frame→True,
PlotLabel→"Basis Functions\n β1(t), β2(t) and Derivatives",
ImageSize→Small]
```

*Out[23]=*

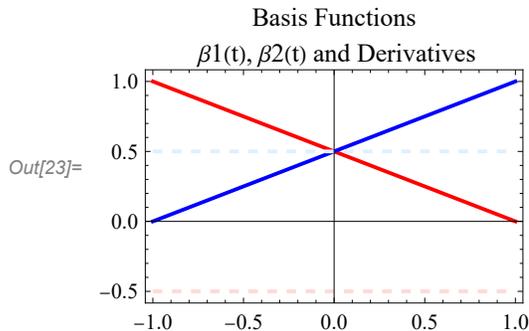

▲ **Figure 2.** The piecewise linear basis functions, $\beta_1(t)$ and $\beta_2(t)$, and their derivatives on a typical boundary element.

### ■ *Normal Vector and Fundamental Solutions*

As required for the BEM, we will define the unit outward normal vector on the boundary, denoted by *NormalVe*c, the fundamental solution, and the fundamental flux. We previously used the symbol $w$ for the fundamental flux in Eq. (2) so that $\frac{\partial w}{\partial n}$ would denote the fundamental flux. Now the notation *uFundamental* and *fluxFundamental*, respectively, will be used. The fundamental solution will be used to approximate the unknown boundary outward normal flux and, subsequently, the values of the primary unknown in the domain.

```mathematica
In[24]:= NormalVec[i_] = {y[i+1]-y[i],x[i]-x[i+1]} / Sqrt[(y[i+1]-y[i])^2+(x[i]-x[i+1])^2];

        uFundamental[xVar_,yVar_,ξ_,η_] = -Log[Sqrt[(xVar-ξ)^2+(yVar-η)^2]] / (2*π);

        fluxFundamental[xVar_,yVar_,ξ_,η_] = -{(xVar-ξ),(yVar-η)} / (2*π*((xVar-ξ)^2+(yVar-η)^2));
```

We assign the boundary node coordinates, $(x(i), y(i)) = BEMNodeCoords(i)$. Recall that the $i^{th}$ boundary element is delineated by the consecutive boundary nodes indexed $i$ and $i + 1$. Then we can define *uExact*, the exact solution at each the boundary node. Recall that the exact boundary solution is known since this is a Dirichlet problem.

```mathematica
In[27]:= Table[
          {x[i],y[i]} = BEMNodeCoords[i] =
            {Cos[2*π*i/nBoundaryNodes],Sin[2*π*i/nBoundaryNodes]}//N,
          {i,1,nBoundaryNodes}];
```

```mathematica
In[28]:= BEMNodeCoords[nBoundaryNodes+1] = BEMNodeCoords[1];

        u[nBoundaryNodes+1] = u[1];

        uExact[BEMNodeCoords[nBoundaryNodes+1]] =
            uExact[BEMNodeCoords[1]];
```

```mathematica
In[31]:= Table[
          u[i] = uExact[BEMNodeCoords[i]],
          {i,1,nBoundaryNodes+1}];
```

We view the exact solution and the boundary nodes.

```
In[32]:=  plotBoundaryPoints =
          ListPointPlot3D[
          Table[{BEMNodeCoords[i][[1]],BEMNodeCoords[i][[2]],u[i]},
          {i,1,nBoundaryNodes}], ImageSize→Small,
          PlotRange→All,PlotStyle→{Black,PointSize[0.03]},
          AxesLabel→{"x","y","u"}, PlotLabel→"Boundary Nodes"];

          Show[plotCircExact,plotBoundaryPoints,
          PlotLabel→ "Exact Analytical Solution\nwith Boundary Nodes",
          AxesLabel→{"x","y","u"}, ImageSize→Small]
```

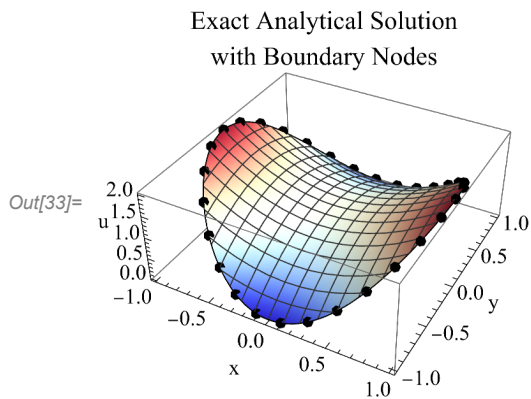

Out[33]=

▲ **Figure 3.** BEM boundary nodes overlaid on the exact solution on the unit disk domain.

The construction of the numerical solution over the entire domain, the disk, will be based only on data at the boundary nodes of the domain, a circle in this case.

## ◻ BEM: Solving for the Outward Normal Flux on the Boundary

### ▪ *BEM Formulation for the Outward Normal Boundary Flux*

The first of the two main loops in our BEM template includes the calculation of the outward normal boundary flux at the boundary. Below we construct this first main loop with the equations to be solved at each of the boundary nodes, as in Eq. (4). This becomes a system of equations with one equation for each boundary collocation node, $(\xi, \eta)$.

For this loop, the numerical integration strategies, rules, and preprocessors become crucial. In order to increase the efficiency in computing time and accuracy, we must choose proper solvers for numerical integration.

The potentially singular integrals, i.e., those integrals for which the field points $(x, y)$ and the collocation point $(\xi, \eta)$ are in the same boundary element, are zero due to orthogonality of the boundary element, a line segment, with respect to the normal vector on that element. This is the reason that, in the following do loop, we use an *If* command to set such integration values equal to zero.

```
In[34]:=  flux[nBoundaryNodes+1] = flux[1];

          Table[flux[i],{i,1,nBoundaryNodes+1}];
```

```
In[36]:=  Do[
          {ξ,η}={x[k],y[k]};

          Do[
          If[((k≠i)&&(k≠i+1))&&({k,i}≠{1,nBoundaryNodes}),
          left[i,k] =
           NIntegrate[
              (fluxFundamental[xVar[i,t],yVar[i,t],ξ,η].
                 NormalVec[i]) * (u[i]*β1[t]+u[i+1]*β2[t]) *
                   Sqrt[(D[xVar[i,t],t])²+(D[yVar[i,t],t])²],
                      {t,-1,1}, Method→integralMethod],
          left[i,k] = 0];

          right[i,k] =
           (flux[i] * NIntegrate[
              uFundamental[xVar[i,t],yVar[i,t],ξ,η] *
                 β1[t] * Sqrt[(D[xVar[i,t],t])²+(D[yVar[i,t],t])²],
                    {t,-1,1}, Method→integralMethod] ) +
           (flux[i+1] * NIntegrate[
              uFundamental[xVar[i,t],yVar[i,t],ξ,η] *
                 β2[t] * Sqrt[(D[xVar[i,t],t])²+(D[yVar[i,t],t])²],
                    {t,-1,1}, Method→integralMethod]),
          {i,1,nBoundaryNodes}],

          {k,1,nBoundaryNodes}]
```

### ▪ Solving for the Outward Normal Boundary Flux

Once the loop for calculating the outward normal boundary flux is constructed, we can rewrite our system as a BEM equation, using the form specified by Eq. (4).

```
In[37]:=  Do[
          BEMequation[k] =
          (FreeTermCoefficient * u[k] + ∑_{i=1}^{nBoundaryNodes} left[i,k])
                == (∑_{i=1}^{nBoundaryNodes} right[i,k]),
          {k,1,nBoundaryNodes}]
```

In this notebook we use an equation-based formulation to solve the BEM system of equations, but using an equivalent matrix-vector formulation with *LinearSolve* solves the system more efficiently.

Now that we have our system in a correct form for the BEM, we can numerically solve for the outward normal boundary flux at each BEM boundary node, and assign it to the quantity *fluxBEM*.

```
In[38]:=  solForFlux =
          Solve[
          Table[BEMequation[i],{i,1,nBoundaryNodes}],
            Table[flux[i],{i,1,nBoundaryNodes}]]//Flatten;

          Table[fluxBEM[k] = solForFlux[[k,2]],{k,1,nBoundaryNodes}];
          fluxBEM[nBoundaryNodes+1]=fluxBEM[1];
```

```
In[41]:=  Table[
          qExact[BEMNodeCoords[k]],{k,1,nBoundaryNodes}];
```

The absolute and relative outward normal boundary flux errors can be computed and displayed, if desired.

```
In[42]:=  Table[
          Abs[fluxBEM[k] -
            qExact[BEMNodeCoords[k]]],{k,1,nBoundaryNodes}];
```

## ◻ BEM: Solving for the Primary Unknown in the Interior

### ■ Grid for Interior Nodes

Here we will we use the user-specified number of interior nodes in the x-direction and y-direction to create a square grid circumscribing the domain. This is one possible and simple way just one possible way to do this. In the square circumscribing the disk domain, the number of interior nodes in the x-direction is constrained to equal the number of interior nodes in the y-direction, previously specified by the user as *tempInteriorNodes*. The total number of these nodes is denoted by *allNodes*.

Recall that an interior domain discretization is not required in the BEM, and this is a major advantage of this method. However, we compute an interior grid here for convenience, so that we will have sufficient uniformly located interior nodes to observe the accuracy of this BEM implementation.

```
In[43]:=   nxnode = nynode = tempInteriorNodes;
           allNodes = nxnode*nynode;
```

The step-sizes in both the x and y-directions, $h_x$ and $h_y$, respectively, are the same. These are hard-coded in the following commands,

$$h_x = h_y = \frac{2}{\text{nxnode} - 1}.$$

```
In[45]:=   {hx,hy}={ 2/(nxnode-1) , 2/(nynode-1) };
```

To determine the number of interior nodes that lie within the unit disk, as desired, we perform a loop to obtain the node coordinates of these interior nodes and store the total number of nodes of interest as *nInteriorNodes*.

```
In[46]:=  kInteriorNode=0;

          Do[
          Do[
          allCoords[i,j] = {i*hx-1, j*hy-1};

          k = i + j*nynode;

          If[(allCoords[i,j][[1]])^2 + (allCoords[i,j][[2]])^2 < 1,
          (* x^2 + y^2 < 1 ,
          inequality sign since boundary node
          coordinates already computed *)

          interiorNodeIndex[k] = ++kInteriorNode;

          InteriorCoords[kInteriorNode] = allCoords[i,j]];
          nInteriorNodes = kInteriorNode,
          {i,0,nxnode-1}],
          {j,0,nynode-1}];
```

Next we check the total number of nodes that are actually contained within this unit disk.

```
In[48]:=  nInteriorNodes
```

Out[48]= 69

To view the interior nodes, we start by creating a unit circle plot, where we will have nodes on the boundary, nodes in the interior, and some nodes that will be ignored outside of the unit disk. To ensure that only the nodes in the interior are selected, we will label and color them differently. A plot overlaying these three items is shown for 2D visualization; some nodes outside of the domain boundary are shown but are not used.

```
In[49]:=  circleplot =
          ParametricPlot[{Cos[z],Sin[z]}, {z,-π,π},
          PlotLabel→"Boundary Node Grid",
          PlotStyle→{Blue}, Frame→True];

          gridplot = ListPlot[
          Partition[
          Table[
          Table[
          allCoords[i,j],
          {i,0,nxnode-1}],
          {j,0,nynode-1}]//Flatten,
          2],PlotMarkers→"●",PlotStyle→Black];

          ActualNodeLabels =
          Table[Text[Style["*",Blue,14],InteriorCoords[k]+{0.0,0.0}],
             {k,1,nInteriorNodes}];

          BoundaryNodeLabels=
          Table[Text[Style["*",Red,18],BEMNodeCoords[k]],
             {k,1,nBoundaryNodes}];

          Show[circleplot,gridplot,
             Graphics[ActualNodeLabels],Graphics[BoundaryNodeLabels],
              PlotLabel→"Boundary and Domain Nodes",
              ImageSize→300]
```

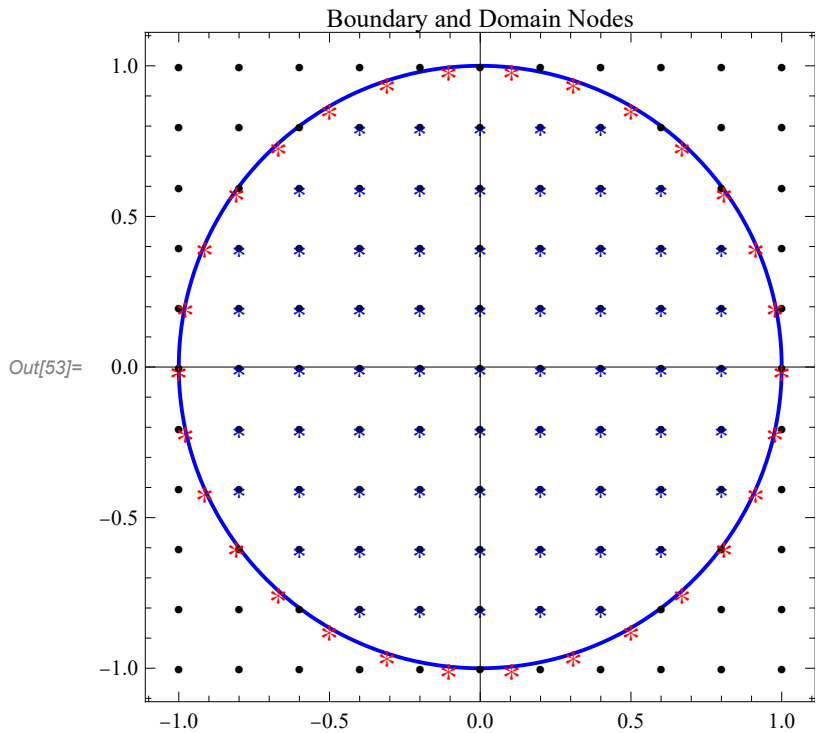

▲ **Figure 4.** BEM nodes on the domain boundary and in interior; black nodes are unused exterior nodes, red stars are boundary nodes, and blue stars are interior nodes.

Using a rectangular geometry and generating the above rectangular grid or mesh may simplify the generation of interior nodes. However, we emphasize this is done solely for convenience here. Discretization of the domain is not required in the BEM. Only the boundary needs to be discretized and meshed. This fact is a major advantage in the use of the BEM, as opposed to finite element and finite difference methods, both of which require discretization, meshing, and solving a system of equations that applies on the entire domain.

### ▪ BEM Formulation for u in the Interior Domain

Similarly to the creation of the loop for solving for the outward normal boundary flux, we now create the second main loop in order to solve for the BEM solution in the domain interior. The distinction is that now we have collocation points, $(\xi, \eta)$, in the interior, as in Eq. (5).

As in the first main loop, numerical integration in this loop can be time-consuming and inefficient. The use of strategies, rules, and preprocessors becomes important for the efficient execution of our BEM technique.

There are no singular integrals since the collocation points are now in the domain interior.

```
In[54]:=  Do[
          {ξ,η} = InteriorCoords[k];
          
          Do[
          left[i,k] = NIntegrate[
             (fluxFundamental[xVar[i,t],yVar[i,t],ξ,η].NormalVec[i]) *
                (u[i]*β1[t]+u[i+1]*β2[t]) *
                   Sqrt[(D[xVar[i,t],t])^2+(D[yVar[i,t],t])^2],
             {t,-1,1}, Method→integralMethod];
          
          right[i,k]=
          fluxBEM[i] * NIntegrate[
             uFundamental[xVar[i,t],yVar[i,t],ξ,η] *
                β1[t] * Sqrt[(D[xVar[i,t],t])^2+(D[yVar[i,t],t])^2],
                   {t,-1,1}, Method→integralMethod] +
          fluxBEM[i+1] * NIntegrate[
             uFundamental[xVar[i,t],yVar[i,t],ξ,η] *
                β2[t] * Sqrt[(D[xVar[i,t],t])^2+(D[yVar[i,t],t])^2],
                   {t,-1,1}, Method→integralMethod],
          {i,1,nBoundaryNodes}],
          
          {k,1,nInteriorNodes}]
```

- ### Constructing the Solution for u in the Interior Domain from the Boundary Solution

We can now create a system of equations to be solved for the BEM solution in the interior of the disk. From Eqs. (5) and (6), we create a discrete set of values that will constitute our numerical solution for *u* in the domain interior, *uBEMInterior*.

```
In[55]:=  Do[
          uBEMInterior[k] = ∑_{i=1}^{nBoundaryNodes} right[i,k] - ∑_{i=1}^{nBoundaryNodes} left[i,k],
          {k,1,nInteriorNodes}]
```

# ◻ BEM: Numerical Results

### ■ *Plotting the BEM Solution*

We now have the BEM flux at the boundary nodes and the BEM solution at the selected nodes in the domain interior. With these results, we can now create graphics in order to graphically compare our BEM solution with the exact solution. First we define the interior data, and secondly we define the boundary data.

```
In[56]:=   interiorData =
           Table[{InteriorCoords[k][[1]],InteriorCoords[k][[2]],
           uBEMInterior[k]},
           {k,1,nInteriorNodes}];

           boundaryData =
           Table[
           {x[i]//N,
           y[i]//N,
           u[i]//N},
           {i,1,nBoundaryNodes}];

           allData = Join[boundaryData,interiorData]//N;
```

In order to partially check that the joined table is correct, we compare the lengths of two tables. The lengths should agree.

```
In[59]:=   totalNodes = nInteriorNodes+nBoundaryNodes;
           totalNodes == (allData//Length)
```

Out[60]= True

To compare the exact solution and the BEM solution, we will plot the discrete set of BEM values and the overlay of the BEM solution on the exact solution.

```mathematica
In[61]:= allBEMpointsplot = 
  ListPointPlot3D[allData,
   PlotRange->All, ColorFunction->"Rainbow",
   AxesLabel->{"x","y","u"}, PlotLabel->"BEM\nSolution",
   RegionFunction->Function[{x,y,z},x^2+y^2<=1],
   PlotStyle->{PointSize->0.02},ImageSize->Small];

allBEMplot = 
  ListPlot3D[allData,
   PlotRange->All, ColorFunction->"TemperatureMap",
   AxesLabel->{"x","y","u"}, PlotLabel->"Smoothed\nBEM Solution",
   RegionFunction->Function[{x,y,z},x^2+y^2<=1], ImageSize->Small];

overlayplot = 
  Show[plotCircExact,allBEMplot,allBEMpointsplot,
   AxesLabel->{"x","y","u"},ImageSize->Small,
   PlotLabel->"Exact and BEM Solutions"];

GraphicsGrid[{{allBEMpointsplot,allBEMplot},{overlayplot}},
 ImageSize->Full]
```

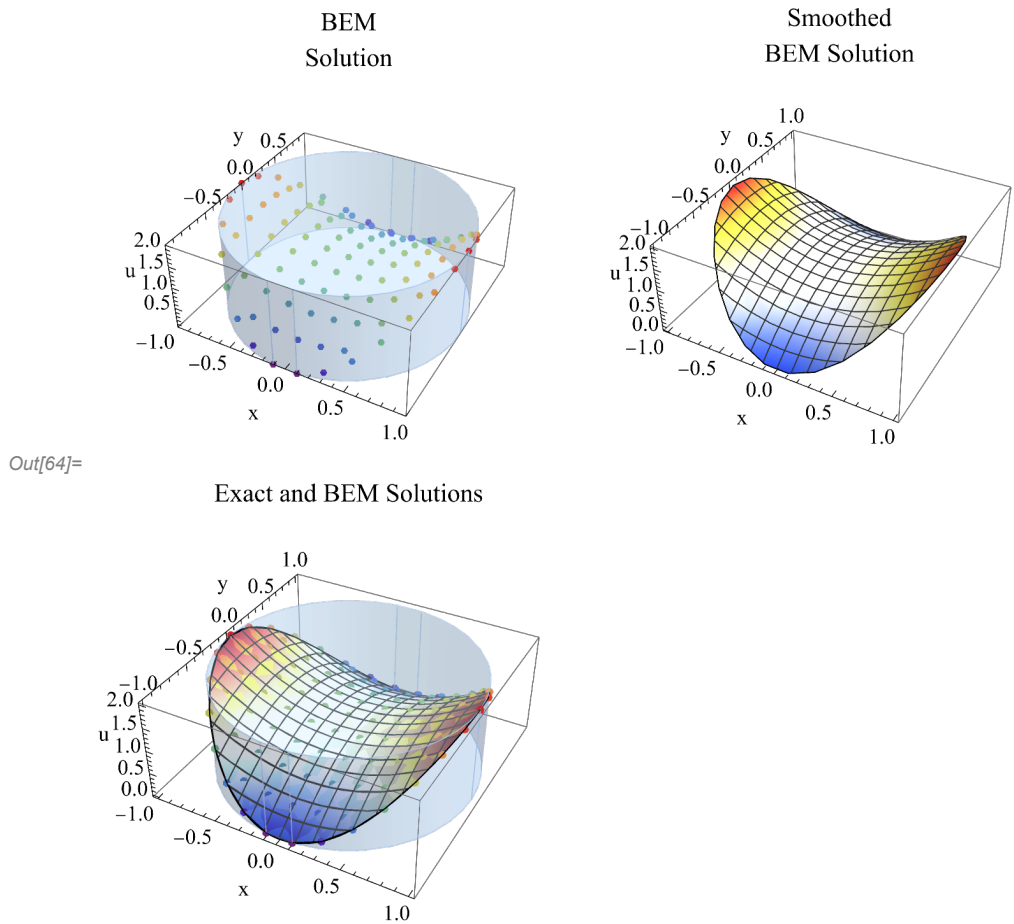

*Out[64]=*

▲ **Figure 5.** (Top Left) A plot of the discrete BEM solution. (Top Right) A plot of the smoothed BEM solution. (Bottom Left) An overlay of the two plots.

### ■ *Numerical Results*

We also create a table for numerical confirmation of the accuracy of the BEM solution.

*Boundary Element Method for the Dirichlet Problem for Laplace's Equation on a Disk*    **23**```
In[65]:=  tableSolResults=
          Labeled[TableForm[
          Table[
          {k,
          InteriorCoords[k][[1]],
          InteriorCoords[k][[2]],
          uBEMInterior[k],
          uExact[InteriorCoords[k]]//N},
          {k,1,nInteriorNodes}],

          TableAlignments→Center, TableSpacing→Automatic,
          TableHeadings→{None,
             {"Interior\nNode #","x","y",
              "BEM\nSolution","Exact\nSolution"}}]//Chop,
          Framed[
          "BEM Solution vs. Exact Solution: Numerical Results"],Top]
```

| BEM Solution vs. Exact Solution: Numerical Results | | | | |
|---|---|---|---|---|
| Interior Node # | x | y | BEM Solution | Exact Solution |
| 1 | $-\frac{2}{5}$ | $-\frac{4}{5}$ | 0.522255 | 0.52 |
| 2 | $-\frac{1}{5}$ | $-\frac{4}{5}$ | 0.402681 | 0.4 |
| 3 | 0 | $-\frac{4}{5}$ | 0.362854 | 0.36 |
| 4 | $\frac{1}{5}$ | $-\frac{4}{5}$ | 0.402681 | 0.4 |
| 5 | $\frac{2}{5}$ | $-\frac{4}{5}$ | 0.522255 | 0.52 |
| 6 | $-\frac{3}{5}$ | $-\frac{3}{5}$ | 0.99997 | 1. |
| 7 | $-\frac{2}{5}$ | $-\frac{3}{5}$ | 0.80089 | 0.8 |
| 8 | $-\frac{1}{5}$ | $-\frac{3}{5}$ | 0.681424 | 0.68 |
| 9 | 0 | $-\frac{3}{5}$ | 0.641602 | 0.64 |
| 10 | $\frac{1}{5}$ | $-\frac{3}{5}$ | 0.681424 | 0.68 |
| 11 | $\frac{2}{5}$ | $-\frac{3}{5}$ | 0.80089 | 0.8 |
| 12 | $\frac{3}{5}$ | $-\frac{3}{5}$ | 0.99997 | 1. |
| 13 | $-\frac{4}{5}$ | $-\frac{2}{5}$ | 1.47798 | 1.48 |
| 14 | $-\frac{3}{5}$ | $-\frac{2}{5}$ | 1.19911 | 1.2 |
| 15 | $-\frac{2}{5}$ | $-\frac{2}{5}$ | 1. | 1. |
| 16 | $-\frac{1}{5}$ | $-\frac{2}{5}$ | 0.880534 | 0.88 |

Printed by Wolfram Mathematica Student Edition



*Out[65]=*

| | | | | |
|---|---|---|---|---|
| 16 | $\frac{?}{5}$ | $\frac{?}{5}$ | ? | ? |
| 17 | 0 | $-\frac{2}{5}$ | 0.840712 | 0.84 |
| 18 | $\frac{1}{5}$ | $-\frac{2}{5}$ | 0.880534 | 0.88 |
| 19 | $\frac{2}{5}$ | $-\frac{2}{5}$ | 1. | 1. |
| 20 | $\frac{3}{5}$ | $-\frac{2}{5}$ | 1.19911 | 1.2 |
| 21 | $\frac{4}{5}$ | $-\frac{2}{5}$ | 1.47798 | 1.48 |
| 22 | $-\frac{4}{5}$ | $-\frac{1}{5}$ | 1.59734 | 1.6 |
| 23 | $-\frac{3}{5}$ | $-\frac{1}{5}$ | 1.31858 | 1.32 |
| 24 | $-\frac{2}{5}$ | $-\frac{1}{5}$ | 1.11947 | 1.12 |
| 25 | $-\frac{1}{5}$ | $-\frac{1}{5}$ | 1. | 1. |
| 26 | 0 | $-\frac{1}{5}$ | 0.960178 | 0.96 |
| 27 | $\frac{1}{5}$ | $-\frac{1}{5}$ | 1. | 1. |
| 28 | $\frac{2}{5}$ | $-\frac{1}{5}$ | 1.11947 | 1.12 |
| 29 | $\frac{3}{5}$ | $-\frac{1}{5}$ | 1.31858 | 1.32 |
| 30 | $\frac{4}{5}$ | $-\frac{1}{5}$ | 1.59734 | 1.6 |
| 31 | $-\frac{4}{5}$ | 0 | 1.63716 | 1.64 |
| 32 | $-\frac{3}{5}$ | 0 | 1.3584 | 1.36 |
| 33 | $-\frac{2}{5}$ | 0 | 1.15929 | 1.16 |
| 34 | $-\frac{1}{5}$ | 0 | 1.03982 | 1.04 |
| 35 | 0 | 0 | 1. | 1. |
| 36 | $\frac{1}{5}$ | 0 | 1.03982 | 1.04 |
| 37 | $\frac{2}{5}$ | 0 | 1.15929 | 1.16 |
| 38 | $\frac{3}{5}$ | 0 | 1.3584 | 1.36 |
| 39 | $\frac{4}{5}$ | 0 | 1.63716 | 1.64 |
| 40 | $-\frac{4}{5}$ | $\frac{1}{5}$ | 1.59734 | 1.6 |
| 41 | $-\frac{3}{5}$ | $\frac{1}{5}$ | 1.31858 | 1.32 |
| 42 | $-\frac{2}{5}$ | $\frac{1}{5}$ | 1.11947 | 1.12 |
| 43 | $-\frac{1}{5}$ | $\frac{1}{5}$ | 1. | 1. |
| 44 | 0 | $\frac{1}{5}$ | 0.960178 | 0.96 |
| 45 | $\frac{1}{5}$ | $\frac{1}{5}$ | 1. | 1. |
| 46 | $\frac{2}{5}$ | $\frac{1}{5}$ | 1.11947 | 1.12 |
| 47 | $\frac{3}{5}$ | $\frac{1}{5}$ | 1.31858 | 1.32 |
| 48 | $\frac{4}{5}$ | $\frac{1}{5}$ | 1.59734 | 1.6 |





| | | | | |
|---|---|---|---|---|
| 49 | $-\frac{4}{5}$ | $\frac{2}{5}$ | 1.47798 | 1.48 |
| 50 | $-\frac{3}{5}$ | $\frac{2}{5}$ | 1.19911 | 1.2 |
| 51 | $-\frac{2}{5}$ | $\frac{2}{5}$ | 1. | 1. |
| 52 | $-\frac{1}{5}$ | $\frac{2}{5}$ | 0.880534 | 0.88 |
| 53 | 0 | $\frac{2}{5}$ | 0.840712 | 0.84 |
| 54 | $\frac{1}{5}$ | $\frac{2}{5}$ | 0.880534 | 0.88 |
| 55 | $\frac{2}{5}$ | $\frac{2}{5}$ | 1. | 1. |
| 56 | $\frac{3}{5}$ | $\frac{2}{5}$ | 1.19911 | 1.2 |
| 57 | $\frac{4}{5}$ | $\frac{2}{5}$ | 1.47798 | 1.48 |
| 58 | $-\frac{3}{5}$ | $\frac{3}{5}$ | 0.99997 | 1. |
| 59 | $-\frac{2}{5}$ | $\frac{3}{5}$ | 0.80089 | 0.8 |
| 60 | $-\frac{1}{5}$ | $\frac{3}{5}$ | 0.681424 | 0.68 |
| 61 | 0 | $\frac{3}{5}$ | 0.641602 | 0.64 |
| 62 | $\frac{1}{5}$ | $\frac{3}{5}$ | 0.681424 | 0.68 |
| 63 | $\frac{2}{5}$ | $\frac{3}{5}$ | 0.80089 | 0.8 |
| 64 | $\frac{3}{5}$ | $\frac{3}{5}$ | 0.99997 | 1. |
| 65 | $-\frac{2}{5}$ | $\frac{4}{5}$ | 0.522255 | 0.52 |
| 66 | $-\frac{1}{5}$ | $\frac{4}{5}$ | 0.402681 | 0.4 |
| 67 | 0 | $\frac{4}{5}$ | 0.362854 | 0.36 |
| 68 | $\frac{1}{5}$ | $\frac{4}{5}$ | 0.402681 | 0.4 |
| 69 | $\frac{2}{5}$ | $\frac{4}{5}$ | 0.522255 | 0.52 |

### ■ Error Analysis

In order to gain more insight into the accuracy of our BEM solution, we calculate the absolute error and the magnitude of the relative error between the BEM solution and the exact solution. The absolute error is calculated as

$$\text{Absolute Error} = |u_{\text{BEM}} - u_{\text{Exact}}|,$$

and the magnitude of the relative error is calculated as

$$\text{Relative Error Magnitude} = \left| \frac{u_{\text{BEM}} - u_{\text{Exact}}}{u_{\text{Exact}}} \right|.$$





```
In[66]:= tableSolError=
Labeled[TableForm[
Table[
{k,
InteriorCoords[k][[1]],
InteriorCoords[k][[2]],

Abs[uBEMInterior[k]-
  uExact[InteriorCoords[k]]],

Abs[(uBEMInterior[k]-
  uExact[InteriorCoords[k]]) /
    uExact[InteriorCoords[k]]]
},
{k,1,nInteriorNodes}]/.∞→Nothing,

TableAlignments→Center, TableSpacing→Automatic,
TableHeadings→{None,
{"Interior\nNode #","x","y",
"Absolute\nError","Relative\nError Magnitude"}}]//Chop,
Framed["BEM Solution vs. Exact Solution: Error Results"],Top]
```

| BEM Solution vs. Exact Solution: Error Results | | | | |
|---|---|---|---|---|
| Interior Node # | x | y | Absolute Error | Relative Error Magnitude |
| 1 | $-\frac{2}{5}$ | $-\frac{4}{5}$ | 0.00225491 | 0.00433637 |
| 2 | $-\frac{1}{5}$ | $-\frac{4}{5}$ | 0.00268107 | 0.00670268 |
| 3 | 0 | $-\frac{4}{5}$ | 0.00285358 | 0.00792662 |
| 4 | $\frac{1}{5}$ | $-\frac{4}{5}$ | 0.00268107 | 0.00670268 |
| 5 | $\frac{2}{5}$ | $-\frac{4}{5}$ | 0.00225491 | 0.00433637 |
| 6 | $-\frac{3}{5}$ | $-\frac{3}{5}$ | 0.0000297287 | 0.0000297287 |
| 7 | $-\frac{2}{5}$ | $-\frac{3}{5}$ | 0.000889741 | 0.00111218 |
| 8 | $-\frac{1}{5}$ | $-\frac{3}{5}$ | 0.00142394 | 0.00209403 |
| 9 | 0 | $-\frac{3}{5}$ | 0.00160194 | 0.00250304 |
| 10 | $\frac{1}{5}$ | $-\frac{3}{5}$ | 0.00142394 | 0.00209403 |
| 11 | $\frac{2}{5}$ | $-\frac{3}{5}$ | 0.000889741 | 0.00111218 |
| 12 | $\frac{3}{5}$ | $-\frac{3}{5}$ | 0.0000297287 | 0.0000297287 |





Out[66]=

| | | | | |
|---|---|---|---|---|
| 13 | $-\frac{4}{5}$ | $-\frac{2}{5}$ | 0.00201729 | 0.00136303 |
| 14 | $-\frac{3}{5}$ | $-\frac{2}{5}$ | 0.000890193 | 0.000741828 |
| 15 | $-\frac{2}{5}$ | $-\frac{2}{5}$ | $6.52661 \times 10^{-10}$ | $6.52661 \times 10^{-10}$ |
| 16 | $-\frac{1}{5}$ | $-\frac{2}{5}$ | 0.000533981 | 0.000606797 |
| 17 | 0 | $-\frac{2}{5}$ | 0.000711974 | 0.000847588 |
| 18 | $\frac{1}{5}$ | $-\frac{2}{5}$ | 0.000533981 | 0.000606797 |
| 19 | $\frac{2}{5}$ | $-\frac{2}{5}$ | $6.52661 \times 10^{-10}$ | $6.52661 \times 10^{-10}$ |
| 20 | $\frac{3}{5}$ | $-\frac{2}{5}$ | 0.000890193 | 0.000741828 |
| 21 | $\frac{4}{5}$ | $-\frac{2}{5}$ | 0.00201729 | 0.00136303 |
| 22 | $-\frac{4}{5}$ | $-\frac{1}{5}$ | 0.00265873 | 0.00166171 |
| 23 | $-\frac{3}{5}$ | $-\frac{1}{5}$ | 0.00142395 | 0.00107875 |
| 24 | $-\frac{2}{5}$ | $-\frac{1}{5}$ | 0.00053398 | 0.000476768 |
| 25 | $-\frac{1}{5}$ | $-\frac{1}{5}$ | $2.48745 \times 10^{-10}$ | $2.48745 \times 10^{-10}$ |
| 26 | 0 | $-\frac{1}{5}$ | 0.000177994 | 0.000185411 |
| 27 | $\frac{1}{5}$ | $-\frac{1}{5}$ | $2.48745 \times 10^{-10}$ | $2.48745 \times 10^{-10}$ |
| 28 | $\frac{2}{5}$ | $-\frac{1}{5}$ | 0.00053398 | 0.000476768 |
| 29 | $\frac{3}{5}$ | $-\frac{1}{5}$ | 0.00142395 | 0.00107875 |
| 30 | $\frac{4}{5}$ | $-\frac{1}{5}$ | 0.00265873 | 0.00166171 |
| 31 | $-\frac{4}{5}$ | 0 | 0.0028422 | 0.00173305 |
| 32 | $-\frac{3}{5}$ | 0 | 0.00160194 | 0.0011779 |
| 33 | $-\frac{2}{5}$ | 0 | 0.000711973 | 0.00061377 |
| 34 | $-\frac{1}{5}$ | 0 | 0.000177993 | 0.000171148 |
| 35 | 0 | 0 | $3.4768 \times 10^{-10}$ | $3.4768 \times 10^{-10}$ |
| 36 | $\frac{1}{5}$ | 0 | 0.000177993 | 0.000171148 |
| 37 | $\frac{2}{5}$ | 0 | 0.000711973 | 0.00061377 |
| 38 | $\frac{3}{5}$ | 0 | 0.00160194 | 0.0011779 |
| 39 | $\frac{4}{5}$ | 0 | 0.0028422 | 0.00173305 |
| 40 | $-\frac{4}{5}$ | $\frac{1}{5}$ | 0.00265873 | 0.00166171 |
| 41 | $-\frac{3}{5}$ | $\frac{1}{5}$ | 0.00142395 | 0.00107875 |
| 42 | $-\frac{2}{5}$ | $\frac{1}{5}$ | 0.00053398 | 0.000476768 |
| 43 | $-\frac{1}{5}$ | $\frac{1}{5}$ | $2.48745 \times 10^{-10}$ | $2.48745 \times 10^{-10}$ |
| 44 | 0 | $\frac{1}{5}$ | 0.000177994 | 0.000185411 |
| 45 | $\frac{1}{5}$ | $\frac{1}{5}$ | ... | ... |





| 45 | $\frac{1}{5}$ | $\frac{1}{5}$ | $2.48745 \times 10^{-10}$ | $2.48745 \times 10^{-10}$ |
|---|---|---|---|---|
| 46 | $\frac{2}{5}$ | $\frac{1}{5}$ | 0.00053398 | 0.000476768 |
| 47 | $\frac{3}{5}$ | $\frac{1}{5}$ | 0.00142395 | 0.00107875 |
| 48 | $\frac{4}{5}$ | $\frac{1}{5}$ | 0.00265873 | 0.00166171 |
| 49 | $-\frac{4}{5}$ | $\frac{2}{5}$ | 0.00201729 | 0.00136303 |
| 50 | $-\frac{3}{5}$ | $\frac{2}{5}$ | 0.000890193 | 0.000741828 |
| 51 | $-\frac{2}{5}$ | $\frac{2}{5}$ | $6.52661 \times 10^{-10}$ | $6.52661 \times 10^{-10}$ |
| 52 | $-\frac{1}{5}$ | $\frac{2}{5}$ | 0.000533981 | 0.000606797 |
| 53 | 0 | $\frac{2}{5}$ | 0.000711974 | 0.000847588 |
| 54 | $\frac{1}{5}$ | $\frac{2}{5}$ | 0.000533981 | 0.000606797 |
| 55 | $\frac{2}{5}$ | $\frac{2}{5}$ | $6.52661 \times 10^{-10}$ | $6.52661 \times 10^{-10}$ |
| 56 | $\frac{3}{5}$ | $\frac{2}{5}$ | 0.000890193 | 0.000741828 |
| 57 | $\frac{4}{5}$ | $\frac{2}{5}$ | 0.00201729 | 0.00136303 |
| 58 | $-\frac{3}{5}$ | $\frac{3}{5}$ | 0.0000297287 | 0.0000297287 |
| 59 | $-\frac{2}{5}$ | $\frac{3}{5}$ | 0.000889741 | 0.00111218 |
| 60 | $-\frac{1}{5}$ | $\frac{3}{5}$ | 0.00142394 | 0.00209403 |
| 61 | 0 | $\frac{3}{5}$ | 0.00160194 | 0.00250304 |
| 62 | $\frac{1}{5}$ | $\frac{3}{5}$ | 0.00142394 | 0.00209403 |
| 63 | $\frac{2}{5}$ | $\frac{3}{5}$ | 0.000889741 | 0.00111218 |
| 64 | $\frac{3}{5}$ | $\frac{3}{5}$ | 0.0000297287 | 0.0000297287 |
| 65 | $-\frac{2}{5}$ | $\frac{4}{5}$ | 0.00225491 | 0.00433637 |
| 66 | $-\frac{1}{5}$ | $\frac{4}{5}$ | 0.00268107 | 0.00670268 |
| 67 | 0 | $\frac{4}{5}$ | 0.00285358 | 0.00792662 |
| 68 | $\frac{1}{5}$ | $\frac{4}{5}$ | 0.00268107 | 0.00670268 |
| 69 | $\frac{2}{5}$ | $\frac{4}{5}$ | 0.00225491 | 0.00433637 |





```
In[67]:= absoluteError =
        Table[
         {InteriorCoords[k][[1]],
          InteriorCoords[k][[2]],
          Abs[uBEMInterior[k] -
             uExact[InteriorCoords[k]]]
         },{k,1,nInteriorNodes}];

        relativeError =
        Table[
         {InteriorCoords[k][[1]],
          InteriorCoords[k][[2]],
          Abs[(uBEMInterior[k] -
             uExact[InteriorCoords[k]]) /
            uExact[InteriorCoords[k]]]
         },{k,1,nInteriorNodes}]/.∞→Nothing;
```

We calculate the maximum and the mean for both the absolute and relative primary unknown errors. Similar error analysis can be performed for the outward normal boundary flux, if desired.

```
In[69]:= MaxAbs =
          Max[Table[Abs[uBEMInterior[k] - uExact[InteriorCoords[k]]],
              {k,1,nInteriorNodes}]];
        MaxRel =
          Max[Table[Abs[(uBEMInterior[k] -
            uExact[InteriorCoords[k]]) /
              uExact[InteriorCoords[k]]],
                {k,1,nInteriorNodes}]];
        MeanAbs =
          Mean[Table[Abs[uBEMInterior[k] -
            uExact[InteriorCoords[k]]],
              {k,1,nInteriorNodes}]];
        MeanRel =
          Mean[Table[Abs[(uBEMInterior[k] -
            uExact[InteriorCoords[k]]) /
              uExact[InteriorCoords[k]]],
                {k,1,nInteriorNodes}]];
```

The table below shows the maximum absolute error, the maximum magnitude of the relative error, the mean absolute error, and the mean magnitude of the relative error for the pri-



mary unknown.

```
In[73]:=  errorData = Grid[Prepend[List[{MaxAbs,MaxRel,MeanAbs,MeanRel}],
          Style[#,{9,Bold}]&/@{"Maximum\nAbsolute Error",
          "Maximum Relative\nError Magnitude",
          "Mean\nAbsolute Error","Mean Relative\nError Magnitude"}],
          Spacings→{Automatic,2},Frame→All,
          Background→{None,{{LightBlue,None}}}]
```

Out[73]=

| Maximum Absolute Error | Maximum Relative Error Magnitude | Mean Absolute Error | Mean Relative Error Magnitude |
|---|---|---|---|
| 0.00285358 | 0.00792662 | 0.00119869 | 0.00161062 |

▲ **Table 1.** Maximum and mean values for the absolute error and the magnitude of the relative error for the primary unknown.

To view the absolute and relative error, we display simultaneously the absolute error and the magnitude of the relative error. It is important to note the scale of these plots and the ac curacy of the BEM solution compared to the exact solution. As usual, grids can be refined to obtain more accuracy.





```
In[74]:=  absoluteErrorplot =
          ListPlot3D[absoluteError,
          PlotRange→All, ColorFunction→Hue,
          AxesLabel→{"x","y","Error "},
          PlotLabel→"BEM Absolute\nError"];

          relativeErrorplot =
          ListPlot3D[relativeError,
          PlotRange→All, ColorFunction→Hue,
          AxesLabel→{"x","y","Error "},
          PlotLabel→"BEM Relative\nError Magnitude"];

          GraphicsRow[{absoluteErrorplot,relativeErrorplot},
             ImageSize→Full]
```

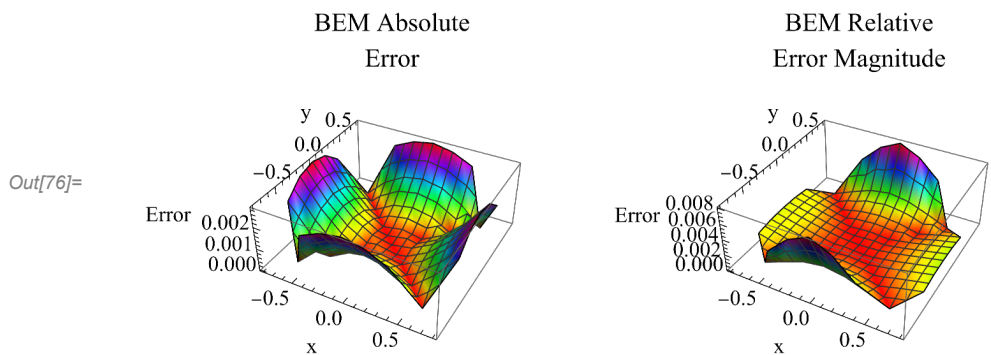

▲ **Figure 6.** (Left) Absolute error for the BEM solution on the domain. (Right) The magnitude of the relative error of the BEM solution on the domain.

By increasing the number of boundary points the absolute error and relative error magnitude are reduced significantly. In general, we see greater error on the boundary compared to the center of the domain, but refining the boundary grid will reduce this.

## ◻ Notebook Timing

To conclude our BEM template, we time the overall runtime of the notebook. This is done by just subtracting the start absolute time from the current absolute time.





```
In[77]:= AbsoluteTime[] - StartTime   (* notebook run time in seconds *)

Out[77]= 13.8982405
```

## ■ Summary

Our main goal was to create a template notebook for using the BEM to solve a Dirichlet problem for Laplace's equation on a disk domain. We started by defining key parameters for the method, including the fundamental quantities and numerical parameters such as the number of boundary nodes. It is also important to select numerical integration strategies, rules, and matrix preprocessors in order to significantly increase the efficiency, both in timing and accuracy. By reformulating the problem as an integral problem and solving for the outward normal boundary flux and interior values, we created a numerical solution to the test problem. Furthermore, graphics were created to confirm the accuracy, and tables for the numerical values and the error values were generated.

## ■ Results

Numerically and graphically we confirmed that the BEM solution is a reasonable approximation for the exact solution to the Dirichlet problem for Laplace's equation on a disk. The piecewise linear basis functions provided reasonable results. For more accuracy, finer grids and/or higher degree basis functions can be used. We computed the absolute error and the magnitude of the relative error, and found reasonable accuracy, considering relatively low degree basis functions and a relatively coarse grid. Of course, what is considered to be "reasonable accuracy" depends on the specific application. In this notebook, this error analysis was done for the primary unknown; this can similarly be done for the flux.

It is important to note that for the sample solution 4, *uExact4*, there is a large discontinuity in the values of the primary unknown as the x-values approach 0. These cause large absolute and relative errors, and degrade accuracy. Nonetheless, the BEM solution is still a reasonable approximation of the exact solution on the rest of the domain.

In terms of timing, there is a drastic increase in notebook run time as we increase the number of nodes, especially when we increase the number of interior nodes. The use of different preprocessors and strategies helps significantly to alleviate the expense of the numerical integration.





# ◼ Conclusion

The BEM is a numerical technique growing in use and applications across science and engineering fields. Its main advantage is that discretization is required only along the boundary of the domain, as opposed to other commonly used methods such as finite differences and finite elements, for which the entire domain must be discretized with a grid. Some BEM disadvantages are that fundamental solutions must be known, and the linear systems to be solved have full coefficient matrices. More theoretical treatment of the BEM mathematics and programming can be found in [12, 13, 14, 15].

Using *LinearSolve* as opposed to *Solve* and formulating the BEM in matrix-vector form, as opposed to equation-based form, will improve efficiency. The numerical integration could be more efficient if it were geared specifically towards the particular application. This code could be modified and adapted for use with a non-uniform boundary grid, other PDEs or differential operators provided the associated fundamental solutions are available, other boundary conditions including Neumann or Robin boundary conditions, non-homogeneous PDEs by implementing the dual reciprocity method for instance, and 3D problems. The Laplace-Transform BEM is an extension of BEM for parabolic and hyperbolic PDEs.

Our template serves as an illustration for solving a Dirichlet problem for Laplace's equation on a disk. The user can specify the number of nodes and the numerical integration strategy in order to improve the accuracy and efficiency of the programming. We are able to approximate the solution using the BEM in an efficient and accurate manner for this 2D model problem.

# ◼ Acknowledgments

I would like to express my gratitude to Dr. Shirley Pomeranz for her guidance and patience and for her continuous support and dedication toward my research and my academic development at The University of Tulsa.

# ◼ Code Availability

The code will be made publicly available in the author's repository:
github.com/misaelmmorales/BEM-Solver





# ◾ References


[1]  W. T. Ang, *A Beginner's Course in Boundary Element Methods*, United States: Universal Publishers, 2007.

[2]  R. Pecher, and J. F. Stanislav, "Boundary element techniques in petroleum reservoir simulation," *Journal of Petroleum Science and Engineering*, **17**, 1997 pp. 353–366. doi.org/10.1016/S0920-4105(96)00066-6

[3]  G. Moridis, and D.L. Reddell, "The Laplace transform boundary element (LTBE) method for the solution of diffusion-type equations " *Boundary elements*, **XIII**, 1991 pp. 83–97. doi.org/10.1007/978-94-011-3696-9_7

[4]  G. Moridis, "Alternative formulations of the Laplace Transform Boundary Element (LTBE) numerical method for the solution of diffusion-type equations " *Boundary element technology*, **VII**, 1992. pp. 815–833. doi.org/10.1007/978-94-011-2872-8_55

[5]  N. M. Al-Ajmi, M. Ahmadi, E. Ozkan, and H. Kazemi. "Numerical Inversion of Laplace Transforms in the Solution of Transient Flow Problems with Discontinuities " *Society of Petroleum Engineers*, 2008. 10.2118/116255-MS

[6]  T. LaForce. "PE281 Boundary Element Method Course Notes" from Stanford University—Energy Resources Engineering department. http://web.stanford.edu/class/energy281.

[7]  C. Constanda, *Solution Techniques for Elementary Partial Differential Equations*, United States: CRC Press, 2018.

[8]  D. L. Young, et. al., "The method of fundamental solutions and condition number analysis for inverse problems of Laplace equation," C*omputers and Mathematics with Applications*, **55** (6), 2008. pp. 1189–1200. doi.org/10.1016/j.camwa.2007.05.015

[9]  A. Iglesia and H. Power, "Boundary Elements with Mathematica," *WIT Transactions on Engineering Sciences*, **15**, 1997. 10.2495/IMS970331

[10] L. Gaul, M. Kögl, and M. Wagner, *Boundary Element Methods for Engineers and Scientists: An Introductory Course with Advanced Topics*, Germany: Springer Berlin-Heidelberg, 2013.

[11] "Wolfram Mathematica Tutorial Collection: Advanced Numerical Integration in Mathematica" from Wolfram *Library Archive*. United States. 2008. library.wolfram.com/infocenter/Books/8504/AdvancedNumericalIntegrationInMathematica.pdf

[12] C. A. Brebbia, and J. Dominguez, *Boundary Elements: An Introductory Course,* Sydney Grammar School Press, 1994.

[13] G. Beer, et. al., *The Boundary Element Method with Programming: For Engineers and Scientists*, Austria: Springer Vienna, 2008.

[14] S. A. Sauter, and C. Schwab, *Boundary Element Methods*, Germany: Springer Berlin-Heidelberg, 2011.

[15] E. P. Stephan, and J. Gwinner, *Advanced Boundary Element Methods: Treatment of Boundary Value, Transmission and Contact Problems*, Germany: Springer International Publishing, 2018.






## About the Authors

Misael M. Morales is an M.S. student in Applied Mathematics at The University of Tulsa. He holds a B.S. in Applied Mathematics and a B.S. in Petroleum Engineering also from The University of Tulsa. His current interests include computational science and mathematics, reservoir modeling and simulation, and data science applications.

Shirley Pomeranz is an Associate Professor of Mathematics at The University of Tulsa. She obtained her PhD in Mathematics from the University of Massachusetts, M.S. in Mathematics from the University of Connecticut, M.S. in Physics from New York University, and B.S. in Physics from Barnard College of Columbia University. She specializes in finite and boundary element methods and numerical analysis.


**Misael M. Morales**
*Department of Mathematics*
*The University of Tulsa*
*800 S. Tucker Drive*
*Tulsa, OK 74104*
*mmm206@utulsa.edu*
† *misaelmorales@utexas.edu*

**Dr. Shirley Pomeranz**
*Department of Mathematics*
*The University of Tulsa*
*800 S. Tucker Drive*
*Tulsa, OK 74104*
*pomeranz@utulsa.edu*